\documentclass[12pt,a4paper]{scrartcl}
\usepackage{amssymb}
\usepackage{amsmath}
\usepackage[ngerman,american]{babel}    %Trennungen, Schriftsatz; Neue deutsche Rechtschreibung
\usepackage[T1]{fontenc}       %Umlaute, Sonderzeichen...
\usepackage[utf8]{inputenc}
\usepackage{graphicx}          %Paket um Grafiken einzubinden. Evtl. muss unter Windows
\usepackage{multicol}          %Paket f�r mehrspaltige Dokumente
\usepackage{float}             %Paket um Bilder in Fliess-Umgebungen einzubinden
\usepackage{color}
\usepackage{amsmath}           %F�r "Formel"-Umgebungen: Wird nachher zum sch�nen Darstellen von 
\usepackage{amssymb}           %Formeln gebrauch (Befehl align)
\usepackage{amsthm}            %F�r Theorem-Umgebungen zu benutzen
\usepackage{listings}          %Paket um Quellcode darzustellen
\usepackage{lastpage}          %Spezielles Paket um die Seiten zu z�hlen
\usepackage{listings}       % Quelltexte einbinden
\usepackage{subfig}
\usepackage{setspace}
\usepackage{mathtools}
\usepackage{todonotes}

\theoremstyle{theorem}

\newtheorem{lemma}{Lemma}[section]			%lemma: counted within the section
\newtheorem{prop}[lemma]{Proposition}			%proposition: counted within the section, same counter as lemma

\newcommand{\V}{\ensuremath{\mathrm{V}}}
\usepackage{url}                               %Hier folgen die Kopfzeilentexte
\usepackage{xspace}

\title{A Sylow theorem for the integral group ring of $\operatorname{PSL}(2,q)$}

\author{Leo Margolis}
\begin{document}
\maketitle
\noindent

{\bf Abstract:} For $G= \operatorname{PSL}(2,p^f)$ denote by $\mathbb{Z}G$ the integral group ring over $G$ and by $V(\mathbb{Z}G)$ the group of units of augmentation 1 in $\mathbb{Z}G.$ Let $r$ be a prime different from $p$. Using the so called HeLP-method we prove that units of $r$-power order in $V(\mathbb{Z}G)$ are rationally conjugate to elements of $G.$ As a consequence we prove that subgroups of prime power order in $V(\mathbb{Z}G)$ are rationally conjugate to subgroups of $G$, if $p=2$ or $f = 1.$\\

Let $G$ be a finite group and $\mathbb{Z}G$ the integral group ring over $G$. Denote by $V(\mathbb{Z}G)$ the group of units of augmentation 1 in $\mathbb{Z}G$, i.e. those units whose coefficients sum up to $1.$ We say that a finite subgroup $U$ of $V(\mathbb{Z}G)$ is rationally conjugate to a subgroup $W$ of $G$, if there exists a unit $x \in \mathbb{Q}G$ such that $x^{-1}Ux=W.$ The question if some, or even all, finite subgroups of $V(\mathbb{Z}G)$ are rationally conjugate to subgroups of $G$ was proposed by H. J. Zassenhaus in the '60s and published in \cite{Zassenhaus}. This so called Zassenhaus Conjectures motivated a lot of research. E.g. A. Weiss proved the strongest version, that all finite subgroups of $V(\mathbb{Z}G)$ are rationally conjugate to subgroups of $G$, provided $G$ is nilpotent \cite{Weiss88} \cite{Weiss91}. K. W. Roggenkamp and L. L. Scott obtained a counterexample \cite{RoggenkampScott} to this strong conjecture. The version, which asks whether all finite cyclic subgroups of $V(\mathbb{Z}G)$ are rationally conjugate to subgroups of $G$, the so called First Zassenhaus Conjecture, is however still open, see e.g. \cite{HertweckEdinb}, \cite{CyclicByAbelian}. Though mostly solvable groups were considered when studying such questions, there are some results available for series of non-solvable groups. E.g. a work on the symmetric groups \cite{Peterson} or for Lie-groups of small rank \cite{Bleher99}. The groups $\operatorname{PSL}(2,q)$, which are also the object of study in this paper, found also some special attention in \cite{Wagner}, \cite{HertweckBrauer}, \cite{HertweckHoefertKimmi}, \cite{AndreasKimmi} or in \cite{PSL2p3}.\\ 
In this paper we will limit our attention to "Sylow-like" results, i.e. to finite $p$-subgroups of $V(\mathbb{Z}G).$ We say that a weak Sylow theorem holds for $V(\mathbb{Z}G)$, if every finite $p$-subgroup of $V(\mathbb{Z}G)$ is isomorphic to some subgroup of $G,$ and that a strong Sylow theorem holds for $V(\mathbb{Z}G)$, if every finite $p$-subgroup of $V(\mathbb{Z}G)$ is even rationally conjugate to a subgroup of $G.$ First Sylow-like results for integral group rings were obtained in \cite{RoggenkampKimmi}. Later M. A. Dokuchaev and S. O. Juriaans proved a strong Sylow theorem for special classes of solvable groups \cite{DokuchaevJuriaans} and M. Hertweck, C. Höfert and W. Kimmerle proved a weak Sylow theorem for $\operatorname{PSL}(2,p^f),$ where $p=2$ or $f \leq 2.$ The results of this article are as follows:\\   

{\bf Theorem 1:} Let $G = \operatorname{PSL}(2,p^f)$, let $r$ be a prime different from $p$ and let $u$ be a torsion unit in $V(\mathbb{Z}G)$ of $r$-power order. Then $u$ is rationally conjugate to a group element.\\

{\bf Theorem 2:} Let $G = \operatorname{PSL}(2,p^f)$ such that $f = 1$ or $p = 2.$ Then any subgroup of prime power order of $V(\mathbb{Z}G)$ is rationally conjugate to a subgroup of $G$, i.e. a strong Sylow theorem holds in $V(\mathbb{Z}G).$\\

\section{HeLP-method and known results}

Let $G$ be a finite group. A very useful notion to study rational conjugacy of torsion units are partial augmentations: Let $u = \sum\limits_{g \in G} a_g g \in \mathbb{Z}G$ and $x^G$ be the conjugacy class of the element $x \in G$ in $G.$ Then $\varepsilon_x(u) = \sum\limits_{g \in x^G} a_g$ is called the {\bf partial augmentation} of $u$ at $x.$ This relates to rational conjugacy via:

\begin{lemma}[{\cite[Th. 2.5]{MarciniakRitterSehgalWeiss}}] Let $u \in \V(\mathbb{Z}G)$ be a torsion unit. Then $u$ is rationally conjugate to a group element if and only if $\varepsilon_x(u^k) \geq 0$ for all $x \in G$ and all powers $u^k$ of $u$. \end{lemma}

It is well known that if $u \neq 1$ is a torsion unit in $\V(\mathbb{Z}G)$, then $\varepsilon_1(u)=0$ by the so called Berman-Higman Theorem \cite[Prop. 1.4]{SehgalBook}. If $\varepsilon_x(u) \neq 0$, then the order of $x$ divides the order of $u$ \cite[Th. 2.7]{MarciniakRitterSehgalWeiss}, \cite[Prop. 3.1]{HertweckColloq}. Moreover the exponent of $G$ and of $\V(\mathbb{Z}G)$ coincide \cite[Cor. 4.1]{CohnLivingstone} and if $U$ is a finite subgroup of $\V(\mathbb{Z}G)$ the order of $U$ divides the order of $G$ \cite{ZK} (or \cite[Lemma 37.3]{SehgalBook}). We will use these facts in the following without further mention.\\

Let $u$ be a torsion unit in $V(\mathbb{Z}G)$ of order $n$ and $\zeta$ a primitive complex $n$-th root of unity and let $K$ be some field, whose characteristic $p$ does not divide $n.$ Let $\xi$ be a (not necessarily primitive) complex $n$-th root of unity and let $D$ be a $K$-representation of $G$ with character $\varphi$. Here $\varphi$ is understood as an ordinary or Brauer character. It was first obtained by Luthar and Passi for $K$ having characteristic 0 \cite{LP89} and later generalized by Hertweck for positive characteristic \cite{HertweckBrauer} that the multipilicity of $\xi$ as an eigenvalue of $D(u)$, which we denote by $\mu(\xi,u,\varphi)$ and which is of course a non-negative integer, may be computed as

$$\mu(\xi,u,\varphi) ={1 \over n} \sum_{\substack{d|n \\ d \neq 1}} {\rm{Tr}}_{\mathbb{Q}(\zeta^d)/\mathbb{Q}}(\varphi(u^d)\xi^{-d}) \ + {1 \over n} \smashoperator[r]{\sum_{\substack{x^G \\ x \ p-{\rm{regular}}}}} \varepsilon_x(u){\rm{Tr}}_{\mathbb{Q}(\zeta)/\mathbb{Q}}(\varphi(x)\xi^{-1}),$$

\noindent
where as usual ${\rm{Tr}}_{\mathbb{Q}(\zeta)/\mathbb{Q}}(x) = \smashoperator[r]{\sum\limits_{\sigma \in {\rm{Gal}}(\mathbb{Q}(\zeta)/\mathbb{Q})}}\sigma(x).$\\
The multiplicity of a complex root of unity as an eigenvalue of a matrix over a field of positive characteristic should be understood, here and in the rest of the paper, in the sense of Brauer.\\
If $u$ is of prime power order $r^k$ for the first sum in the expression above we obtain 
$${1 \over n} \sum_{\substack{d|n \\ d \neq 1}} {\rm{Tr}}_{\mathbb{Q}(\zeta^d)/\mathbb{Q}}(\varphi(u^d)\xi^{-d}) = {1 \over r}\mu(\xi^r,u^r,\varphi).$$
Using these formulas to find possible partial augmentations for torsion units in integral group rings of finite groups is today called HeLP-method.
For a diagonalizable matrix $A$ we will write $A \sim (a_1,...,a_n)$, if the eigenvalues of $A$, with multiplicities, are $a_1,...,a_n.$\\

All subgroups of $G = \operatorname{PSL}(2,p^f)$ were first known to Dickson \cite[Theorem 260]{Dickson}. Let $d = {\rm{gcd}}(2,p-1).$ Up to conjugation there is exactly one cyclic group of order $p$, ${p^f + 1 \over d}$ and ${p^f - 1 \over d}$ respectively in $G$ and every element of $G$ lies in a conjugate of such a group. In particular there is only one conjugacy class of involutions in $G.$ The Sylow $p$-subgroups are elementary-abelian, the Sylow subgroups for all other primes, which are odd, are cyclic and if $p \neq 2$ the Sylow 2-subgroup is dihedral or a Kleinian four-group. There are $d$ conjugacy classes of elements of order $p$. If $g \in G$ is not of order $p$ or 2 its only distinct conjugate in $\langle g \rangle$ is $g^{-1}.$ We denote by $a$ a fixed element of order $p^f-1 \over d$ and by $b$ a fixed element of order $p^f+1 \over d.$\\

The modular representation theory of $\operatorname{PSL}(2,p^f)$ in defining characteristic is well known. All irreducible representations were first given by R. Brauer and C. Nesbitt \cite{BrauerNesbitt}. The explicit Brauer table of $\operatorname{SL}(2,p^f)$, which contains the Brauer table of $\operatorname{PSL}(2,p^f)$, may be found in \cite{SrinivasanPSL}. In particular all characters are real valued since any $p$-regular element is conjugate to its inverse. However, I was not able to find the following Lemma in the literature, except, without proof, in Hertwecks preprint \cite{HertweckBrauer}, so a short proof is included.

\begin{lemma}\label{preps}
Let $G= \operatorname{PSL}(2,p^f)$ and $d={\rm{gcd}}(2,p-1).$ There are $p$-modular representations of $G$ given by $\Theta_0, \Theta_1, \Theta_2,...$ such that there is a primitive $p^f - 1 \over d$-th root of unity $\alpha$ and a primitive $p^f+1 \over d$-th root of unity $\beta$ satisfying
\begin{align*}
\Theta_k(b) &\sim (1,\beta, \beta^{-1}, \beta^2, \beta^{-2},...,\beta^k, \beta^{-k}), \nonumber \\
\Theta_k(a) &\sim (1,\alpha, \alpha^{-1}, \alpha^2, \alpha^{-2},...,\alpha^k, \alpha^{-k}) \nonumber
\end{align*}
for every $k \in \mathbb{N}_0.$
\end{lemma}

{\it Proof:} The group $\operatorname{SL}(2,p^f)$ acts on the vector space spanned by the homogenous polynomials in two commuting variables $x, \ y$ of some fixed degree $e$ extending the natural operation on the 2-dimensional vector space spanned by $x, \ y$, see e.g. \cite[p. 14-16]{Alperin}. Since 
$\left(\begin{smallmatrix} -1 & 0 \\ 0 & -1 \end{smallmatrix}\right) x^iy^j = (-1)^{i+j}x^iy^j$ this action affords a $\operatorname{PSL}(2,p^f)$-representation if and only if $e$ is even and $p$ is odd or $p=2$. So let from now on $e$ be even.\\
Call this representation $\Theta_{e \over 2}.$ Let $\gamma$ be an eigenvalue of an element in $\operatorname{SL}(2,q)$ mapping onto $a$ under the natural projection from $\operatorname{SL}(2,p^f)$ to $\operatorname{PSL}(2,p^f)$. Then $\Theta_{e \over 2}(a)$ has the same eigenvalues as $\Theta_{e \over 2}\left(\left(\begin{smallmatrix} \gamma & 0 \\  0 & \gamma^{-1} \end{smallmatrix}\right)\right).$ Now $\left(\begin{smallmatrix} \gamma & 0 \\  0 & \gamma^{-1} \end{smallmatrix}\right)x^iy^j = \gamma^{i-j}x^iy^j$, so the eigenvalues are $\{\gamma^{i-j} \ | \ 0 \leq i,j \leq e, \ i+j = e\} = \{(\gamma^{2t} \ | \ {-e \over 2} \leq t \leq {e \over 2}\}.$ Thus setting $\alpha = \gamma^2$ proves the first part of the claim.  \\
Now let $\delta$ be an eigenvalue of an element in $\operatorname{SL}(2,p^f)$ mapping onto $b$ under the natural projection from $\operatorname{SL}(2,p^f)$ to $\operatorname{PSL}(2,p^f)$. The action of $\operatorname{SL}(2,p^f)$ may of course be extended to $\operatorname{SL}(2,p^{2f}).$ So $\Theta_{e \over 2}(b)$ has the same eigenvalues as $\Theta_{e \over 2}\left(\left(\begin{smallmatrix} \delta & 0 \\  0 & \delta^{-1} \end{smallmatrix}\right)\right),$ where the matrix may be seen as an element in $\operatorname{SL}(2,p^{2f}).$ Then doing the same calculations as above and setting $\beta = \delta^2$ proves the Lemma.\\ 

{\bf Notation:} The Brauer-character belonging to a representation $\Theta_i$ from the Lemma above will be denoted by $\varphi_i.$\\

Using the HeLP-method R. Wagner \cite{Wagner} and Hertweck \cite{HertweckBrauer} obtained already some results about rational conjugacy of torsion units of prime power order in $\operatorname{PSL}(2,p^f).$ Part of Wagners result was published in \cite{BovdiHoefertKimmi}.

\begin{lemma}\label{Wagner}\cite{Wagner} (also \cite[Prop. 6.1]{HertweckBrauer}) Let $G = \operatorname{PSL}(2,p^f)$ and $f \leq 2.$ If $u$ is a unit of order $p$ in $V(\mathbb{Z}G)$, then $u$ is rationally conjugate to a group element. \end{lemma}

{\it Remark}: The HeLP-method does not suffice to prove rational conjugacy to a group element of a unit of order $p$ in $V(\mathbb{Z}\operatorname{PSL}(2,p^f)),$ if $p$ is odd and $f\geq 3.$ There is also no other method or idea around how one could e.g. obtain whether a unit of order 3 in $V(\mathbb{Z}\operatorname{PSL}(2,27))$ is rationally conjugate to a group element or not.

\begin{lemma}\label{n=1}\cite[Prop. 6.4]{HertweckBrauer} Let $G = \operatorname{PSL}(2,p^f)$ and let $r$ be a prime different from $p.$ If $u$ is a unit of order $r$ in $V(\mathbb{Z}G)$, then $u$ is rationally conjugate to an element of $G.$ \end{lemma}

\begin{lemma}\label{Bovdi}\cite[Prop. 6.5]{HertweckBrauer} Let $G=\operatorname{PSL}(2,p^f)$, let $r$ be a prime different from $p$ and $u$ a torsion unit in $V(\mathbb{Z}G)$ of order $r^n$. Let $m < n$ and denote by $S$ a set of representatives of conjugacy classes of elements of order $r^m$ in $G$. Then $\sum\limits_{x \in S} \varepsilon_x(u) = 0.$\\
If moreover $g$ is an element of order $r^n$ in $G$, then $\mu(1,u,\varphi) = \mu(1,g,\varphi)$ for every $p$-modular Brauer character $\varphi$ of $G$.
 \end{lemma}

If one is interested not only in cyclic groups the following result is very useful. It may be found e.g. in \cite[Lemma 37.6]{SehgalBook} or in \cite[Lemma 4]{Valenti}.

\begin{lemma}\label{Valenti} Let $G$ be a finite group, $U$ a finite subgroup of $V(\mathbb{Z}G)$ and $H$ a subgroup of $G$ isomorphic to $U$. If $\sigma: U \rightarrow H$ is an isomorphism such that $\chi(u) = \chi(\sigma(u))$ for all $u \in U$ and all irreducible complex characters $\chi$ of $G$, then $U$ is rationally conjugate to $H$.  \end{lemma}

\section{Proof of the results}

We will first sum up some elementary number theoretical facts. The notatian $a \equiv b \pod c$ will mean, that $a$ is congruent $b$ modulo $c.$

\begin{prop}\label{Trace} Let $t$ and $s$ be natural numbers such that $s$ divides $t$ and denote by $\zeta_t$ and $\zeta_s$ a primitive complex $t$-th root of unity and $s$-th root of unity respectively. Then  
$${\rm{Tr}}_{\mathbb{Q}(\zeta_t)/\mathbb{Q}}(\zeta_s) = \mu(s){\varphi(t) \over \varphi(s)},$$
where $\mu$ denotes the Möbius function and $\varphi$ Euler's totient function.
So for a prime $r$ and natural numbers $n,m$ with $m \leq n$ we have
$$ {\rm{Tr}}_{\mathbb{Q}(\zeta_{r^n})/\mathbb{Q}}(\zeta_{r^m}) = \left\{ \begin{array}{lll} r^{n-1}(r-1), & m = 0 \\  -r^{n-1}, & m=1 \\ 0, & m > 1 \end{array}\right. $$
Let  moreover $i$ and $j$ be integers prime to $r$, then
$${\rm{Tr}}_{\mathbb{Q}(\zeta_{r^n})/\mathbb{Q}}(\zeta_{r^m}^i\zeta_{r^m}^{-j}) = \left\{ \begin{array}{lll} r^{n-1}(r-1), & i \equiv j \ (r^m) \\  -r^{n-1}, & i \not\equiv j \ (r^m), \ \ i \equiv j \ (r^{m-1}) \\ 0, & i \not\equiv j \ (r^{m-1}) \end{array}\right.$$
\end{prop}

{\it Proof of Proposition \ref{Trace}:} Let $s = p_1^{f_1}\cdot...\cdot p_k^{f_k}$ be the prime factorisation of $s$. For a natural number $l$ let $I(l) = \{ i \in \mathbb{N} \ | \ 1\leq i \leq l, \ {\rm{gcd}}(i,l)=1\}.$ As is well known, ${\rm{Gal}}(\mathbb{Q}(\zeta_t)/\mathbb{Q}) = \{\sigma_i: \zeta_t \mapsto \zeta_t^i \ | \ i \in I(t) \}.$ From this the case $s=1$ follows immediately. Otherwise we have 
$${\rm{Tr}}_{\mathbb{Q}(\zeta_t)/\mathbb{Q}}(\zeta_s) = \sum\limits_{i \in I(t)} \zeta_s^i = {\varphi(t) \over \varphi(s)}\sum\limits_{i \in I(s)} \zeta_s^i = {\varphi(t) \over \varphi(s)} \prod\limits_{j=1}^k \sum\limits_{i \in I(p_j^{f_j})} \zeta_{p_j^{f_j}}^i. $$
Now $\sum\limits_{i \in I(p_j^{f_j})} \zeta_{p_j^{f_j}}^i =  \left\{ \begin{array}{ll} -1, & f_j = 1 \\  0, & f_j > 1 \end{array}\right.$ and this gives the first formula. The other formulas are special cases of this general formula since $\varphi(r^n) = (r-1)(r^{n-1}).$ \\

{\it Proof of Theorem 1:} Let $G = \operatorname{PSL}(2,p^f)$, let $r$ be a prime different from $p$ and let $u$ be a torsion unit in $V(\mathbb{Z}G)$ of order $r^n.$ Let $\zeta$ be an primitive complex $r^n$-th root of unity and set ${\rm{Tr}}_{\mathbb{Q}(\zeta)/\mathbb{Q}} = {\rm{Tr}}.$ If $n = 1$, then by Lemma \ref{n=1} $u$ is rationally conjugate to an element in $G$, so assume $n \geq 2.$ Assume further that by induction $u^r$ is rationally conjugate to an element in $G$.\\
We will proceed by induction on $m$, where $1\leq m \leq n.$ For a fixed $m$ let $l = {r^m -1  \over 2}$ if $r$ is odd and $l = {r^m - 2 \over 2}$ if $r=2.$ Let $\{x_i \ | \ 1\leq i \leq l, \ {\rm{gcd}}(i,r)=1\}$ be a full set of representatives of conjugacy classes of elements of order $r^m$ in $G$ such that $x_1^i = x_i$ (this is possible by the group theoretical properties of $G$ given above) and let $S$ be a set of representatives of conjugacy classes of elements of $G$ of $r$-power order not greater than $r^n$ containing $x_1,...,x_l.$ The proof will be divided in several steps:
\begin{enumerate}
\item[a)] For $m < n$ we will show $\varepsilon_{x_i}(u) = 0$ for every $1\leq i \leq l,$ i.e. the partial augmentations of $u$ at elements of order $r^m$ vanish. This will be proved by an induction on $k \leq m$ showing
\begin{enumerate}
\item[(i)]  $\varepsilon_{x_i}(u) = \varepsilon_{x_j}(u)$ for $i \equiv \pm j \pod {r^{m-k}}.$
\end{enumerate}
If $r$ is even it suffices to prove this for $k =m - 1,$ if $r$ is odd we will prove it for $k =m.$ It then follows from Lemma \ref{Bovdi} that $\varepsilon_{x_i}(u) = 0$ for all $x_i.$
\item[b)] For $m = n$ after reordering the $x_1,...,x_l$ we will show $\varepsilon_{x_1}(u) = 1$ and $\varepsilon_{x_i}(u) = 0$ for $i \geq 2.$ Together with a) this proves the Theorem. This will be achieved by an induction on $k$ to show several facts:
\begin{enumerate}
\item[(i)] $\varepsilon_{x_1}(u) = 1$ and $\varepsilon_{x_i}(u) = 0$ for $i \equiv \pm 1 \pod {r^{n-k}}, i \neq 1.$
\item[(ii)] $\varepsilon_{x_i}(u) = \varepsilon_{x_j}(u)$ for $i \equiv \pm j \pod {r^{n-k}}$ and $i \not \equiv \pm 1 \pod {r^{n-k}}.$ 
\end{enumerate} 
These statements will be proved for $k = n-1.$ If $r$ is even this is already enough. In case $r$ is odd we will moreover prove that $\sum\limits_{i \equiv \alpha (r)}\varepsilon_{x_i}(u) = 0$ for $\alpha \not \equiv \pm 1 \pod r,$ which then also implies the Theorem.
\end{enumerate}

So let $m$ be a natural number such that $m<n.$ If $m=0$ statement a) is the Berman-Higman Theorem and if $r=2$ and $m=1$ it follows from Lemma \ref{Bovdi} and the fact that there is only one conjugacy class of involutions in $G$. So assume we know $\varepsilon_x(u) = 0$ for $\circ(x) < r^m.$  The representations and the corresponding characters of $G$ from Lemma \ref{preps} will be used freely. Statement (i) in a) is certainly true for $k=0,$ so assume $\varepsilon_{x_i}(u) = \varepsilon_{x_j}(u)$ for $i \equiv \pm j \pod {r^{m-k}}$ for some $k$.  Since $u^r$ is rationally conjugate to a group element, there exists a primitive $r^{n-1}$-th root of unity $\zeta_{r^{n-1}}$ such that
$$\Theta_{r^k}(u^r) \sim (1, \zeta_{r^{n-1}}, \zeta_{r^{n-1}}^{-1}, \zeta_{r^{n-1}}^2, \zeta_{r^{n-1}}^{-2},...,\zeta_{r^{n-1}}^{r^k}, \zeta_{r^{n-1}}^{-r^k}).$$
Now all $p$-modular Brauer characters of $G$ are real valued and thus using the last statement in Lemma \ref{Bovdi} we obtain that $\Theta_{r^k}(u) \sim (1,a_1,a_1^{-1},a_2,a_2^{-1},...,a_{r^k},a_{r^k}^{-1})$, where for every $i$ we have $a_i$ a root of unity such that $a_i^{r^{m-k}} \neq 1.$ So for every primitive $r^{m-k}$-th root of unity $\zeta_{r^{m-k}}$ we have $\mu(\zeta_{r^{m-k}},u,\varphi_{r^k}) = 0.$ Let $\zeta_{r^m}$ be a primitive $r^{m}$-th root of unity such that we have $\Theta_{r^k}(x_1) \sim (1,\zeta_{r^m},\zeta_{r^m}^{-1},...,\zeta_{r^m}^{r^k},\zeta_{r^m}^{-r^k})$ and set $\xi = \zeta_{r^m}^{r^k}.$ Let moreover $\alpha$ be a natural number prime to $r$ such that $1 \leq \alpha \leq l$.\\
Thus $\mu(\xi^\alpha,u,\varphi_{r^k}) =0$ and $\varepsilon_x(u) = 0$ for $\circ(x) < r^m$. From here on a sum over $i,$ if not specified, will always mean a sum over all $i$ satisfying $1\leq i \leq l$ and $r \nmid i.$ Then using the HeLP-method we get
\begin{align}
0 &= \mu(\xi^\alpha,u,\varphi_{r^k}) = {1 \over r}\mu(\xi^{\alpha r},u^r,\varphi_{r^k}) + {1 \over r^n} \sum\limits_{x \in S} \varepsilon_x(u){\rm{Tr}}(\varphi_{r^k}(x)\xi^{-\alpha}) \nonumber \\
&={1 \over r}\mu(\xi^{\alpha r},u^r,\varphi_{r^k}) + {1 \over r^n} \smashoperator[r]{\sum_{\substack{x \in S \\ \circ(x) > r^m}}} \varepsilon_{x}(u){\rm{Tr}}(\varphi_{r^k}(x)\xi^{-\alpha}) + {1 \over r^n}\sum_{\substack{i}} \varepsilon_{x_i}(u){\rm{Tr}}(\varphi_{r^k}(x_i)\xi^{-\alpha}) \nonumber \\
&= {1 \over r}\mu(\xi^{\alpha r},u^r,\varphi_{r^k}) + {1 \over r^n}\sum_{\substack{x \in S}} \varepsilon_{x}(u){\rm{Tr}}(\xi^{-\alpha}) + {1 \over r^n}\sum_{\substack{i}} \varepsilon_{x_i}(u){\rm{Tr}}((\xi^i + \xi^{-i})\xi^{-\alpha}) \nonumber \\
\label{A}
&= {1 \over r}\mu(\xi^{\alpha r},u^r,\varphi_{r^k}) + {{\rm{Tr}}(\xi^{-\alpha}) \over r^n} + {1 \over r^n}\sum_{\substack{i}} \varepsilon_{x_i}(u){\rm{Tr}}((\xi^i + \xi^{-i})\xi^{-\alpha}).
\end{align}
In the third line we used that if $\tilde\zeta$ is a root of unity of $r$-power order such that $\tilde\zeta^{r^{m-k}} \neq 1$, then $\tilde\zeta\xi$  has the same order as $\tilde\zeta$ and so ${\rm{Tr}}(\tilde\zeta\xi) = 0$ by Proposition \ref{Trace}. Thus for an element $x \in S$ of order at least $r^m$ we get 
$${\rm{Tr}}(\varphi_{r^k}(x)\xi^{-\alpha}) = \left\{\begin{array}{ll} {\rm{Tr}}(\xi^{-\alpha}), & \circ(x) > r^m \\ {\rm{Tr}}(\xi^{-\alpha} + (\xi^i + \xi^{-i})\xi^{-\alpha}), & x = x_i \end{array}\right.$$
Note that as $i$ is prime to $r$ the congruence $i \equiv \alpha \pod {r^{m-k}}$ implies $-i \not\equiv \alpha  \pod {r^{m-k}}$ for $r^{m-k} \notin \{1,2\}$ and these exceptions don't have to be considered by our assumptions on $m$ and $k$.\\
There are now two cases to consider. First assume $k < m-1$, so $\xi$ is at least of order $r^2$ and thus ${\rm{Tr}}(\xi) = 0.$ Moreover $\mu(\xi^{\alpha r},u^r,\varphi_{r^k})=0$ and using Proposition \ref{Trace} in (\ref{A}) we obtain
\begin{align}
0 &= {1 \over r^n}\sum_{\substack{i}} \varepsilon_{x_i}(u){\rm{Tr}}((\xi^i + \xi^{-i})\xi^{-\alpha})  \nonumber \\
&= {1 \over r^n} \smashoperator[r]{\sum_{\substack{i \equiv \pm \alpha (r^{m-k})}}}\varepsilon_{x_i}(u)(r^{n-1}(r-1)) + {1 \over r^n} \smashoperator[r]{\sum_{\substack{i \equiv \pm \alpha (r^{m-k-1}) \\ i \not\equiv \pm \alpha (r^{m-k})}}}\varepsilon_{x_i}(u)(-r^{n-1})  \nonumber \\
\label{B}
&= \smashoperator[r]{\sum_{\substack{i \equiv \pm \alpha (r^{m-k})}}}\varepsilon_{x_i}(u) - {1 \over r} \smashoperator[r]{\sum_{\substack{i \equiv \pm \alpha (r^{m-k-1})}}}\varepsilon_{x_i}(u).
\end{align}
So
$$r\sum_{\mathclap{\substack{i \equiv \pm \alpha (r^{m-k})}}}\varepsilon_{x_i}(u) =\smashoperator[r]{\sum_{\substack{i \equiv \pm \alpha (r^{m-k-1}) }}}\varepsilon_{x_i}(u).$$
But since by induction $\varepsilon_{x_i}(u) = \varepsilon_{x_j}(u)$ for $i \equiv \pm j \pod {r^{m-k}}$ the summands on the left hand side are all equal and since changing $\alpha$ by $r^{m-k-1}$ does not change the right hand side of the equation we get $\varepsilon_{x_i}(u) = \varepsilon_{x_j}(u)$ for $i \equiv \pm j \pod {r^{m-k-1}}.$\\
Now consider $k=m-1$, then $\xi$ is a primitive $r$-th root of unity and thus ${\rm{Tr}}(\xi) = -r^{n-1}$ and $\mu(\xi^{\alpha r},u^r,\varphi_{r^k}) = \mu(1,u^r,\varphi_{r^k}) = 1.$ So using Proposition \ref{Trace} in (\ref{A}) we get
\begin{align}
0 &= {1 \over r} + {-r^{n-1} \over r^n} + {1 \over r^n}\smashoperator[r]{\sum_{\substack{\pm i \not\equiv \alpha (r)}}} \varepsilon_{x_i}(u)(-2r^{n-1}) + {1 \over r^n}\smashoperator[r]{\sum_{\substack{\pm i \equiv \alpha (r)}}} \varepsilon_{x_i}(u)(r^{n-1}(r-1) - r^{n-1}) \nonumber \\
\label{E}
&= \smashoperator[r]{\sum_{\substack{ \pm i \equiv \alpha (r)}}} \varepsilon_{x_i}(u) - {2 \over r}\sum_{\substack{i}} \varepsilon_{x_i}(u). 
\end{align}
So 
$$r\smashoperator[r]{\sum\limits_{\substack{\pm i \equiv \alpha (r)}}} \varepsilon_{x_i}(u) = 2\smashoperator[r]{\sum\limits_{\substack{i}}} \varepsilon_{x_i}(u).$$
Now by Lemma \ref{Bovdi} the right side of this equation is zero and by induction all summands on the left side are equal. Hence varying $\alpha$ gives $\varepsilon_{x_i}(u) = 0$ for every $x_i$ and part a) is proved. 

So assume $m = n$. As in the computation above we have 
$$\Theta_{r^k}(u^r) \sim (1,\zeta_{r^{n-1}},\zeta_{r^{n-1}}^{-1},...,\zeta_{r^{n-1}}^{r^k},\zeta_{r^{n-1}}^{-r^k})$$
for some primitive $r^{n-1}$-th root of unity $\zeta_{r^{n-1}}$ and $\Theta_{r^k}(u) \sim (1,a_1,a_1^{-1},a_2,a_2^{-1},...,a_{r^k},a_{r^k}^{-1})$, where $a_i$ are roots of unity such that $a_i^{r^{n-k}} \neq 1$ for $1 \leq i \leq r^k-1$ and $a_{r^k}$ is some primitive $r^{n-k}$-th root of unity. Set $\xi = a_{r^k}$ and reorder the $x_1,...,x_l$ such that $\Theta_1(x_1) \sim \Theta_1(u)$, but still $x_1^i = x_i$. Then $x_1^r$ is rationally conjugate to $u^r.$  We will proceed by induction on $k.$\\
Let $\alpha$ be a natural number prime to $r$ with $1\leq \alpha \leq l.$ Using the HeLP-method and $\varepsilon_{x}(u) = 0$ for $\circ(x) < r^n$ we obtain, doing the same calculations as in (\ref{A}):
\begin{eqnarray}
\label{C}
\mu(\xi^\alpha,u,\varphi_{r^k}) = {1 \over r}\mu(\xi^{\alpha r},u^r,\varphi_{r^k}) + {{\rm{Tr}}(\xi^{-\alpha}) \over r^n} + {1 \over r^n}\sum_{\substack{i}} \varepsilon_{x_i}(u){\rm{Tr}}((\xi^i + \xi^{-i})\xi^{-\alpha}). 
\end{eqnarray}
As $u^r$ is rationally conjugate to $x_1^r$ we know that $\xi^{\pm r}$ are eigenvalues of $\Theta_{r^k}(u^r).$ So we get  
$$\mu(\xi^\alpha, u, \varphi_{r^k}) = \left\{ \begin{array}{ll} 1, & \alpha \equiv \pm 1 \ (r^{n-k}) \\ 0, & {\rm{else}}  \end{array}\right. \ \ \ \ {\rm{and}} \ \ \ \  \mu(\xi^{\alpha r},u^r,\varphi_{r^k}) = \left\{ \begin{array}{ll} 1, & \alpha \equiv \pm 1 \ (r^{n-k-1}) \\ 0, & {\rm{else}}   \end{array} \right.$$
There are now several cases to consider: Statement (ii) from b) is clear for $k=0.$ So let $k < n -1.$ For $i \not \equiv \pm 1 \pod {r^{n-k-1}},$ set $\alpha = i.$ Then $\mu(\xi^\alpha,u,\varphi_{r^k}) = \mu(\xi^{r\alpha},u^r.\varphi_{r^k}) = 0$ and since $\xi$ has order at least $r^2,$ also ${\rm{Tr}}(\xi) = 0.$ Thus we can do the same computations as in (\ref{B}) to obtain (ii) for $k +1$. So (ii) holds for $k=n-1.$\\
To obtain the base case for (i) set $k=0.$ Then $\xi$ is at least of order $r^2,$ i.e. ${\rm{Tr}}(\xi) = 0$, and from (\ref{C}) we obtain (similar to the computation in (\ref{B})):
\begin{eqnarray}
\label{E}
1 = {1 \over r} + \varepsilon_{x_1}(u) - {1 \over r}\smashoperator[r]{\sum_{\substack{i \equiv \pm 1 (r^{n-1})}}} \varepsilon_{x_i}(u)
\end{eqnarray}
and
\begin{eqnarray}
\label{F}
0 = {1 \over r} + \varepsilon_{x_\alpha}(u) - {1 \over r}\smashoperator[r]{\sum_{\substack{i \equiv \pm 1 (r^{n-1})}}} \varepsilon_{x_i}(u)
\end{eqnarray}
for $\alpha \equiv \pm 1 \pod {r^{n-1}}$ and $\alpha \neq 1.$ Substracting \eqref{F} from \eqref{E} gives 
\begin{eqnarray}
\label{D}
1 = \varepsilon_{x_1}(u) - \varepsilon_{x_\alpha}(u)
\end{eqnarray}
for every $\alpha \equiv \pm 1 \pod {r^{n-1}}$ and $\alpha \neq 1.$ Let $t = |\{i \in \mathbb{N} | i \leq l, i \equiv \pm 1 \pod {r^{n-1}}\}|.$ Then summing up \eqref{E} with the equations \eqref{F} for all $\alpha \equiv \pm 1 \pod {r^{n-1}}$ gives
$$1 = {t \over r} +  \smashoperator[r]{\sum_{\substack{i \equiv \pm 1 (r^{n-1})}}} \varepsilon_{x_i}(u) - {t \over r}\smashoperator[r]{\sum_{\substack{i \equiv \pm 1 (r^{n-1})}}} \varepsilon_{x_i}(u) = {t \over r} + (1 - {t \over r})\smashoperator[r]{\sum_{\substack{i \equiv \pm 1 (r^{n-1})}}} \varepsilon_{x_i}(u).$$
So $\smashoperator[r]{\sum\limits_{\substack{i \equiv \pm 1 (r^{n-1})}}} \varepsilon_{x_i}(u) = 1$ and the base case of (i) follows from (\ref{D}).\\
So assume $0 \leq k < n-1.$ Then $\xi$ is at least of order $r^2$, i.e. ${\rm{Tr}}(\xi) = 0$. By induction $\smashoperator[r]{\sum\limits_{\substack{i \equiv \pm 1 \pod {r^{n-k}}}}} \varepsilon_{x_i}(u) = 1$ and for $\alpha \equiv \pm 1 \pod {r^{n-k}}$ from (\ref{C}) computing as in (\ref{B}) we obtain
$$1 = {1 \over r} + \smashoperator[r]{\sum_{\substack{i \equiv \pm 1 (r^{n-k})}}} \varepsilon_{x_i}(u) - {1 \over r}\smashoperator[r]{\sum_{\substack{i \equiv \pm 1 (r^{n-k-1})}}} \varepsilon_{x_i}(u) = {1 \over r} + 1 - {1 \over r}\smashoperator[r]{\sum_{\substack{i \equiv \pm 1 (r^{n-k-1})}}} \varepsilon_{x_i}(u).$$
For $\alpha \not \equiv \pm 1 \pod {r^{n-k}}$ and $\alpha \equiv \pm 1 \pod {r^{n-k-1}}$ we obtain the same way
$$0 = {1 \over r} +\smashoperator[r]{\sum_{\substack{i \equiv \pm \alpha (r^{n-k})}}} \varepsilon_{x_i}(u) - {1 \over r}\smashoperator[r]{\sum_{\substack{i \equiv \pm 1 (r^{n-k-1})}}} \varepsilon_{x_i}(u).$$
Thus subtracting the last equation from the one before gives 
$$1 =  1 - \smashoperator[r]{\sum_{\substack{i \equiv \pm \alpha (r^{n-k})}}} \varepsilon_{x_i}(u).$$
The summands on the right hand side are all equal by (ii), so $\varepsilon_{x_\alpha}(u) = 0,$ as claimed. \\
Finally let $r$ be odd, $k = n - 1$ and $\alpha \not \equiv \pm 1 \pod r.$ Then $\xi$ is of order $r$, i.e. ${\rm{Tr}}(\xi) = -r^{n-1},$ and  $\mu(\xi^{\alpha r},u^r,\varphi_{r^k}) = \mu(1,u^r,\varphi_{r^k}) = 3.$ So from (\ref{C}) computing as in (\ref{E}) we obtain
$$ 0 = {3 \over r} + {-r^{n-1} \over r^n} - {2 \over r}\smashoperator[r]{\sum_{\substack{i}}} \varepsilon_{x_i}(u) + \smashoperator[r]{\sum_{\substack{i \equiv \pm \alpha (r)}}} \varepsilon_{x_i}(u) = \smashoperator[r]{\sum_{\substack{i \equiv \pm \alpha (r)}}} \varepsilon_{x_i}(u).$$
As by (ii) all summands in the last sum are equal, we get $\varepsilon_{x_\alpha}(u) = 0$ and the Theorem is finally proved.\\

{\it Proof of Theorem 2:} Let $G = \operatorname{PSL}(2,p^f)$ such that $f = 1$ or $p = 2.$ Assume first that $r$ is an odd prime, which is not $p,$ and $R$ is an $r$-subgroup of $V(\mathbb{Z}G).$ As every $r$-subgroup of $G$ is cyclic so is $R$ by \cite[Theorem A]{HertweckCpCp} and thus $R$ is rationally conjugate to a subgroup of $G$ by Theorem 1. If $p \neq 2$ and $R$ is a 2-subgroup of $V(\mathbb{Z}G)$, then $R$ is either cyclic or dihedral or a Kleinian four group by \cite[Theorem 2.1]{HertweckHoefertKimmi}. If $R$ is cyclic, then it is rationally conjugate to a subgroup of $G$ by Theorem 1. If $R$ is dihedral or a Kleinian four group let $S = \langle s \rangle$ be a maximal cyclic subgroup of $R.$ Then any element of $R$ outside of $S$ is an involution and moreover $s$ is rationally conjugate to an element $g \in G$ by Theorem 1. $R$ is isomorphic to some subgroup $H$ of $G$, such that the maximal cyclic subgroup of $H$ is generated by $g$. As there is only one conjugacy class of invoultions in $G$ every isomorphism $\sigma$ between $R$ and $H$ mapping $s$ to $g$ satisfies $\chi(\sigma(u)) = \chi(u)$ for every irreducible complex character of $G$ and every $u \in R.$ Thus $R$ is rationally conjugate to $H$ by Lemma \ref{Valenti}.\\
If $p=2$ and $P$ is a 2-subgroup of $V(\mathbb{Z}G)$ then all non-trivial elements of $P$ are involutions, so $P$ is elementary abelian and the order of $P$ divides the order of $G$. As there is again only one conjugacy class of involutions in $G$ every isomorphism $\sigma$ between $P$ and a subgroup of $G$ isomorphic with $P$ satisfies $\chi(\sigma(u)) = \chi(u)$ for every irreducible complex character of $G$ and every $u \in P.$ So $P$ is rationally conjugate to a subgroup of $G$ by Lemma \ref{Valenti}. Finally if $p$ is odd and $P$ is a $p$-subgroup of $V(\mathbb{Z}G)$, then $P$ is cyclic of order $p$ and thus rationally conjugate to a subgroup of $G$ by Lemma \ref{Wagner}.\\

{\bf Remark:} Let $G= \operatorname{PSL}(2,p^f)$ and let $n$ be a number prime to $p.$ The structure of the Brauer table of $G$ in defining characteristic yields immidiately, that if we can prove that a unit $u \in V(\mathbb{Z}G)$ of order $n$ is rationally conjugate to an element in $G$ applying the HeLP-method to the Brauer table, then these calculations will hold over any $\operatorname{PSL}(2,q)$, if $n$ and $q$ are coprime. In this sense it would be interesting, and seems actually achievable, to determine a subset $A_{p^f}$ of $\mathbb{N}$ such that we can say: The HeLP-method proves that a unit $u \in V(\mathbb{Z}G)$ of order $n$ is rationally conjugate to an element in $G$ if and only if $n \in A_{p^f}.$ Test computations yield the conjecture that $A_{p^f}$  actually contains all odd numbers prime to $p$. If this turned out to be true this would yield, using the results in \cite{HertweckBrauer}, the First Zassenhaus Conjecture for the groups $\operatorname{PSL}(2,p),$ where $p$ is a Fermat- or Mersenne prime.\\
Other interesting questions concerning torsion units of the integral group ring of $ G = \operatorname{PSL}(2,p^f)$ were mentioned at the end of \cite{HertweckHoefertKimmi} and are still open today: If the order of $u \in V(\mathbb{Z}G)$ is divisible by $p$, is $u$ of order $p$? Are units of order $p$ rationally conjugate to elements of $G$? Are $p$-subgroups in $V(\mathbb{Z}G)$ necessarily abelian? For the last question a positive answer in case $f=3$ is given in \cite{PSL2p3}.\\

{\bf Acknowledgement:} The computations given above were all done by hand, but some motivating computations were done using a GAP-implementation of the HeLP-algorithm written by the author in collaboration with Andreas Bächle \cite{HeLP}.\\
I also thank the referee for many suggestions improving the readability of the paper.

\bibliographystyle{amsalpha}
\bibliography{giganty}

\providecommand{\bysame}{\leavevmode\hbox to3em{\hrulefill}\thinspace}
\providecommand{\MR}{\relax\ifhmode\unskip\space\fi MR }
\providecommand{\MRhref}[2]{%
  \href{http://www.ams.org/mathscinet-getitem?mr=#1}{#2}
}
\providecommand{\href}[2]{#2}
\begin{thebibliography}{MRSW87}

\bibitem[Alp86]{Alperin}
J.~L. Alperin, \emph{Local representation theory}, Cambridge Studies in
  Advanced Mathematics, vol.~11, Cambridge University Press, Cambridge, 1986,
  Modular representations as an introduction to the local representation theory
  of finite groups.

\bibitem[BHK04]{BovdiHoefertKimmi}
V.~Bovdi, C.~H{\"o}fert, and W.~Kimmerle, \emph{On the first {Z}assenhaus
  conjecture for integral group rings}, Publ. Math. Debrecen \textbf{65}
  (2004), no.~3-4, 291--303.

\bibitem[BK11]{AndreasKimmi}
A.~B{\"a}chle and W.~Kimmerle, \emph{On torsion subgroups in integral group
  rings of finite groups}, J. Algebra \textbf{326} (2011), 34--46.

\bibitem[Ble99]{Bleher99}
Frauke~M. Bleher, \emph{Finite groups of {L}ie type of small rank}, Pacific J.
  Math. \textbf{187} (1999), no.~2, 215--239.

\bibitem[BM15a]{HeLP}
A.~B{\"a}chle and L.~Margolis, \emph{{HeLP} -- {H}ertweck-{L}uthar-{P}assi
  method, \textsf{GAP} package, {V}ersion 1.0},
  \url{http://homepages.vub.ac.be/abachle/help/}, 2015.

\bibitem[BM15b]{PSL2p3}
Andreas B{\"a}chle and Leo Margolis, \emph{Torsion subgroups in the units of
  the integral group ring of {{PSL}$(2,p^3)$}}, Arch. Math. (Basel)
  \textbf{105} (2015), no.~1, 1--11.

\bibitem[BN41]{BrauerNesbitt}
R.~Brauer and C.~Nesbitt, \emph{On the modular characters of groups}, Ann. of
  Math. (2) \textbf{42} (1941), 556--590.

\bibitem[CL65]{CohnLivingstone}
James~A. Cohn and Donald Livingstone, \emph{On the structure of group algebras.
  {I}}, Canad. J. Math. \textbf{17} (1965), 583--593.

\bibitem[CMdR13]{CyclicByAbelian}
Mauricio Caicedo, Leo Margolis, and {\'A}ngel del R{\'{\i}}o, \emph{Zassenhaus
  conjecture for cyclic-by-abelian groups}, J. Lond. Math. Soc. (2) \textbf{88}
  (2013), no.~1, 65--78.

\bibitem[Dic01]{Dickson}
Leonard~Eugene Dickson, \emph{Linear groups: {W}ith an exposition of the
  {G}alois field theory}, Teubner, Leipzig, 1901.

\bibitem[DJ96]{DokuchaevJuriaans}
Michael~A. Dokuchaev and Stanley~O. Juriaans, \emph{Finite subgroups in
  integral group rings}, Canad. J. Math. \textbf{48} (1996), no.~6, 1170--1179.

\bibitem[Her06]{HertweckColloq}
Martin Hertweck, \emph{On the torsion units of some integral group rings},
  Algebra Colloq. \textbf{13} (2006), no.~2, 329--348.

\bibitem[Her07]{HertweckBrauer}
\bysame, \emph{Partial augmentations and {B}rauer character values of torsion
  units in group rings}, arXiv:math.RA/0612429v2, 2004 - 2007.

\bibitem[Her08a]{HertweckEdinb}
\bysame, \emph{Torsion units in integral group rings of certain metabelian
  groups}, Proc. Edinb. Math. Soc. (2) \textbf{51} (2008), no.~2, 363--385.

\bibitem[Her08b]{HertweckCpCp}
\bysame, \emph{Unit groups of integral finite group rings with no noncyclic
  abelian finite {$p$}-subgroups}, Comm. Algebra \textbf{36} (2008), no.~9,
  3224--3229.

\bibitem[HHK09]{HertweckHoefertKimmi}
Martin Hertweck, Christian~R. H{\"o}fert, and Wolfgang Kimmerle, \emph{Finite
  groups of units and their composition factors in the integral group rings of
  the group {${\rm PSL}(2,q)$}}, J. Group Theory \textbf{12} (2009), no.~6,
  873--882.

\bibitem[KR93]{RoggenkampKimmi}
W.~Kimmerle and K.~W. Roggenkamp, \emph{A {S}ylow-like theorem for integral
  group rings of finite solvable groups}, Arch. Math. (Basel) \textbf{60}
  (1993), no.~1, 1--6.

\bibitem[LP89]{LP89}
I.~S. Luthar and I.~B.~S. Passi, \emph{Zassenhaus conjecture for {$A_5$}},
  Proc. Indian Acad. Sci. Math. Sci. \textbf{99} (1989), no.~1, 1--5.

\bibitem[MRSW87]{MarciniakRitterSehgalWeiss}
Z.~Marciniak, J.~Ritter, S.~K. Sehgal, and A.~Weiss, \emph{Torsion units in
  integral group rings of some metabelian groups. {II}}, J. Number Theory
  \textbf{25} (1987), no.~3, 340--352.

\bibitem[Pet76]{Peterson}
Gary~L. Peterson, \emph{Automorphisms of the integral group ring of {$S_{n}$}},
  Proc. Amer. Math. Soc. \textbf{59} (1976), no.~1, 14--18.

\bibitem[Rog91]{RoggenkampScott}
Klaus~W. Roggenkamp, \emph{Observations on a conjecture of {H}ans
  {Z}assenhaus}, Groups---{S}t. {A}ndrews 1989, {V}ol. 2, London Math. Soc.
  Lecture Note Ser., vol. 160, Cambridge Univ. Press, Cambridge, 1991,
  pp.~427--444.

\bibitem[Seh93]{SehgalBook}
S.~K. Sehgal, \emph{Units in integral group rings}, Pitman Monographs and
  Surveys in Pure and Applied Mathematics, vol.~69, Longman Scientific \&
  Technical, Harlow; copublished in the United States with John Wiley \& Sons,
  Inc., New York, 1993, With an appendix by Al Weiss.

\bibitem[Sri64]{SrinivasanPSL}
Bhama Srinivasan, \emph{On the modular characters of the special linear group
  {$SL(2,\,p^{n})$}}, Proc. London Math. Soc. (3) \textbf{14} (1964).

\bibitem[Val94]{Valenti}
Angela Valenti, \emph{Torsion units in integral group rings}, Proc. Amer. Math.
  Soc. \textbf{120} (1994), no.~1, 1--4.

\bibitem[Wag95]{Wagner}
R.~Wagner, \emph{Zassenhausvermutung über die {G}ruppen
  {$\operatorname{PSL}(2,p)$}}, Master's thesis, Universität Stuttgart, Mai
  1995.

\bibitem[Wei88]{Weiss88}
Alfred Weiss, \emph{Rigidity of {$p$}-adic {$p$}-torsion}, Ann. of Math. (2)
  \textbf{127} (1988), no.~2, 317--332.

\bibitem[Wei91]{Weiss91}
\bysame, \emph{Torsion units in integral group rings}, J. Reine Angew. Math.
  \textbf{415} (1991), 175--187.

\bibitem[Zas74]{Zassenhaus}
Hans Zassenhaus, \emph{On the torsion units of finite group rings}, Studies in
  mathematics (in honor of {A}. {A}lmeida {C}osta) ({P}ortuguese), Instituto de
  Alta Cultura, Lisbon, 1974, pp.~119--126.

\bibitem[{\v{Z}}K67]{ZK}
{\`E}.~M. {\v{Z}}mud and G.~{\v{C}}. Kurenno{\u\i}, \emph{The finite groups of
  units of an integral group ring}, Vestnik Har'kov. Gos. Univ. \textbf{1967}
  (1967), no.~26, 20--26.

\end{thebibliography}

\bigskip

Leo Margolis, Fachbereich Mathematik, Universit\"{a}t Stuttgart, Pfaffenwaldring 57, 70569 Stuttgart, Germany.
\emph{leo.margolis@mathematik.uni-stuttgart.de}

\end{document}